\documentclass[twocolumn,a4]{autart}

\usepackage{amssymb, amsfonts, amsmath}
\usepackage{verbatim}

\usepackage{subfigure}
\usepackage{slashbox}
\usepackage[dvipdfm]{graphicx}

\newcommand{\gv}[1]{\ensuremath{\mbox{\boldmath$ #1 $}}}
\newcommand{\grad}[1]{\gv{\nabla} #1}
\newcommand{\R}{\mathbb{R}}
\newcommand{\N}{\mathbb{N}}
\newcommand{\bmat}[1]{\begin{bmatrix}#1\end{bmatrix}}
\newcommand{\aligneq}[1]{
\begin{align*}
#1
\end{align*}
}
\newcommand{\lrr}[1]{\left ( #1 \right )}
\newcommand{\srr}[1]{\left [ #1 \right ]}
\newcommand{\rmatrix}[2]{\lrr{\begin{array}{#1}
#2
\end{array}}}

\newcommand{\Reals}[0]{\mathbb{R}}
\newcommand{\Realsn}[1]{\mathbb{R}^{#1}}

\newtheorem{mydef}{Definition}

\newtheorem{mylem}{Lemma}

\newtheorem{myass}{Assumption}

\newcommand{\theoremnC}[2]{
\vspace{-.1in}
\begin{thm}
{\bf (#1)}  #2
\end{thm}
\vspace{-.1in}
}

\newcommand{\theoremC}[1]{
\vspace{-.1in}
\begin{thm}
#1
\end{thm}
\vspace{-.1in}
}

\newcommand{\ProofC}[1]{

\begin{pf}
#1
\end{pf}

}

\newcommand{\defnc}[1]{
\vspace{-.1in}
\begin{mydef}
#1
\end{mydef}}

\newcommand{\defnn}[2]{
\vspace{-.1in}
\begin{mydef}
{\bf (#1)}  #2
\end{mydef}}

\newcommand{\norm}[1]{\left\| #1 \right\|}

\newcounter{exmpl_counter}
\newcommand{\exmpl}[1]{\par\addtocounter{exmpl_counter}{1}\noindent \textbf{Example \arabic{exmpl_counter}.}  {#1} \par}

\newcommand{\assn}[1]{\begin{myass}
#1
\end{myass}}

\begin{document}

\begin{frontmatter}

\title{{\bf Using SOS for Analysis of Zeno Stability in Hybrid systems with Nonlinearity and Uncertainy}}
\author[iit]{Chaitanya Murti}
\ead{cmurti@hawk.iit.edu}
\author[asu]{Matthew Peet}
\ead{mpeet@asu.edu}

\address[iit]{Illinois Institute of Technology, Chicago, IL}

\address[asu]{Arizona State University, Tempe, AZ}

\begin{abstract}
Hybrid systems exhibit phenomena which do not occur in systems with continuous vector fields. One such phenomenon - Zeno executions - is characterized by an infinite number of discrete events or transitions occurring over a finite interval of time. This phenomenon is not necessarily undesirable and may indeed be used to capture physical phenomena. In this paper, we examine the problem of proving the existence and stability of zero executions. Our approach is to develop a polynomial-time algorithm - based on the sum-of-squares methodology - for verifying the stability of a Zeno execution. We begin by stating Lyapunov-like theorems for local Zeno stability based on existing results. Then, for hybrid systems with polynomial vector fields, we use polynomial Lyapunov functions and semialgebraic geometry (Positivstellensatz results) to reduce the local Lyapunov-like conditions to a convex feasibility problem in polynomial variables. The feasibility problem is then tested using an algorithm for sum-of-squares programming - SOSTOOLS. We also extend these results to hybrid system with parametric uncertainty, where the uncertain parameters lie in a semialgebriac set. We also provide several examples illustrating the use of our technique.
\end{abstract}

\end{frontmatter}
\vspace{-.1in}
\section{Introduction}
\vspace{-.1in}
\noindent  Hybrid systems exhibit both continuous dynamics and discrete or logical transitions, and are  used to model a variety of physical and artificial systems. Examples of systems modeled using hybrid vector fields include electrical systems with switching~\cite{CTA521}, communication networks with queueing~\cite{HespanhaHybCom}, networked control systems~\cite{NCS_SOS}, embedded systems~\cite{HybEmbed}, biological systems~\cite{lygeros_stoch_hybsys}, and air traffic control~\cite{TomlinAir}. \\
\indent Much of the research on hybrid systems involves extending tools for analysis and control of non-hybrid systems to their hybrid counterparts. Examples of this include stability analysis~\cite{BranStability}, observability and controllability~\cite{BemporadObsCont}, and controller synthesis~\cite{Bemp_PWA_synth}. Existence, uniqueness and continuity of solutions for hybrid systems have also been widely studied, e.g.~\cite{lygeros2003dynamical}. Of particular relevance to this paper is the use of Lyapunov-type conditions for stability (e.g.~\cite{ShevPaden}). A common approach to the use of Lyapunov functions for analysis of hybrid system involves  discontinuous or piecewise-continuous Lyapunov functions~\cite{BranMultLyap}. Methods for the construction of piecewise-quadratic functions using Linear Matrix Inequalities can be found in~\cite{JhnsQuadLyap} for systems with a piecewise-affine vector field. Lyapunov methods for robust stability analysis also exist (e.g.~\cite{PettLennStab}). A result on the use of sum-of-squares for stability analysis of system with piecewise-polynomial vector fields can be found in~\cite{PrajPapach} and~\cite{PapachPrajRobust}.
\\
\indent A Zeno execution is a solution of a hybrid model which predicts infinite transitions between discrete states in a finite interval of time. By definition, a Zeno execution (or arc) is a solution to a hybrid system which converges to what is called a Zeno equilibria - a fixed solution which consists of a sort of infinite loop and which has only discrete transitions with no continuous evolution. These executions may arise due to poor modeling or they may represent real physical processes (e.g. a bouncing which comes to rest in finite time).  Properties of hybrid systems with Zeno executions is described in detail in, for example,~\cite{ZhangZunJohKarLygSasZeno}. For control, the importance of understanding Zeno behavior was demonstrated in~\cite{Fuller}, wherein it was shown that the optimal controller for a relay system would necessarily undergo infinitely many transitions in finite time.  \\
\indent While Zeno behavior is not an intrinsically undesirable property, it can have unwanted repercussions. For example, Zeno executions are notoriously difficult to simulate~\cite{ZhangZunJohKarLygSasZeno}. One approach to dealing with this problem is regularization of the system, as discussed in, say,~\cite{JohanssonReg} and~\cite{LifeAftZenoAmes}.
In this paper, we are interested in the problem of predicting whether a Zeno execution will occur. This problem has been studied for some classes of hybrid system. For example, in~\cite{AmesDFQ}, the authors proposed sufficient conditions for existence of Zeno executions in first quadrant hybrid systems. Sufficient conditions for Zeno behavior in hybrid systems with nonlinear vector fields were proposed in \cite{AmesNonlinconst} using a locally flat approximation of the vector field.  Key to the prediction of Zeno executions is the existence of a well-developed Lyapunov theory. For example, a converse Lyapunov result for systems with an isolated Zeno equilibrium was given in~\cite{AmesLyap} (An equivalent result was given in~\cite{TeelGoebLyap}). Non-isolated Zeno equilibria were treated in \cite{AmesLampTAC}.  These results show that stability of a Zeno execution is equivalent to the existence of a Lyapunov function which proves this property.\\
\indent  In this paper, which expands upon our result given in \cite{Chai_ECC2013}, we provide a computationally tractable test for Zeno stability in hybrid systems with semi-algebraic guard sets; piecewise-polynomial vector fields; and piecewise-polynomial transition maps. Specifically, we develop a polynomial-time algorithm for the construction of the  Lyapunov-like functions proposed in~\cite{AmesLyap} and~\cite{TeelGoebLyap} where the functions are piece-wise polynomial of arbitrary fixed degree. We also extend this method to the verification of Zeno stability for systems with parametric uncertainty in the vector fields, guard sets, domains, and reset maps.   \\
\indent The outline of the paper is as follows: in Section 2, we discuss background, including Zeno stability, Lyapunov theory and relevant concepts from optimization and semialgebraic geometry - including sum-of-squares. In Section 3, we use Sum-of-Squares to find a convex approach for construction of Lyapunov functions for Zeno stability. In Section 4, we provide numerical examples. In Section 5, we extend our results to systems with parametric uncertainties and give additional examples.
\vspace{-.1in}
\section{Background}
\vspace{-.1in}
In this section, we provide the following background material. In Subsection~\ref{subsection:sos}, we define Sum-of-Squares polynomials; in Subsection~\ref{subsection:psatz} we introduce definitions and results from real algebraic geometry, including a Positivstellensatz; in Subsection~\ref{subsection:hybsys}, we define a class of hybrid system along with a definition of execution - the solution of a hybrid system; and in Subsection~\ref{subsection:zeno}, we define for Zeno executions, Zeno Equilibria and Zeno stability - as well as a Lyapunov theorem for the latter property.
\vspace{-.1in}
\subsection{Sum of Squares Polynomials}\label{subsection:sos}
\vspace{-.1in}
In this paper, we will be searching for a piecewise-polynomial Lyapunov function to prove stability of a Zeno equilibrium. While convex, this problem is difficult due to the difficulty in parameterizing the set of positive polynomial Lyapunov functions. Indeed, the question of determining whether a polynomial is positive is known to be NP-hard~\cite{blum1988theory}. However, in this paper, we restrict ourselves to a subset of positive polynomial functions known as the \emph{Sum-of-Squares} polynomials. This choice is not conservative in that it has been shown that Sum-of-Squares polynomial Lyapunov functions are necessary and sufficient for stability of nonlinear systems with polynomials vector field~\cite{Peet_TAC2012}.

Let $\mathbf{R}[x]$ denote the ring of polynomials in variables $x = (x_1,..,x_n) \in \R^n$.

\defnn{Sum of Squares Polynomial}{A polynomial $p(x):\Realsn{n}\rightarrow \Reals$ is said to be Sum of Squares (SOS) if there exist polynomials $f_i(x):\Realsn{n}\rightarrow \Reals $ such that
\[p(x) = \sum_i (f_i(x))^2\]
We use $\Sigma_x\subset \mathbf{R}[x]$ to denote the convex cone of polynomials which are SOS.
}
While any SOS polynomial is obviously positive semidefinite, not all positive semidefinite polynomials are Sum-of-Squares. However, as seen in Theorem~\ref{thm:thmsos} (below), we have an efficient test to determine whether a polynomial is SOS.
\noindent\theoremC{\label{thm:thmsos}
For a polynomial, $p$ of degree $2d$, $p\in \Sigma_x$ if and only if there exists a positive semidefinite matrix $Q$, such that
\[p(x) = Z(x)^TQZ(x)\]
where $Z(x)$ is the vector of monomials of degree $d$ or less
}
For a proof, refer to, say, \cite{ParriloThs}.\\
Theorem~\ref{thm:thmsos} shows us that checking whether a polynomial is SOS is equivalent to checking the existence of a positive-semidefinite matrix $Q$ under affine constraints and can therefore be represented as a Linear Matrix Inequality. Thus, while checking polynomial positivity is NP-hard, checking whether a polynomial is SOS is decidable in polynomial time.
\vspace{-.1in}
\subsection{The Positivstellensatz}\label{subsection:psatz}
\vspace{-.1in}
A Positivstellensatz is a result which shows that Sum-of-Squares polynomials can be used to parameterize the cone of polynomials which are positive on a semialgebraic set. For this section, let $\mathbb{N}^0 = \mathbb{N}\cup\{0\}$.
\defnn{Semialgebraic Set}{A semialgebraic set is a set of the form
\aligneq{
&S:=\left \{x\in\Realsn{n}: f_i(x) \geq 0,i=1,...,n_1,\right.\\
&\qquad\qquad\qquad\qquad\qquad \left. h_i(x) = 0,i=1,...,n_3 \right \}
}
where each $f_i\in\mathbf{R}[x]$, and $h_i\in\mathbf{R}[x]$.
}
In simple terms, a semialgebraic set is a set defined by polynomial equalities and inequalities. We now look at the set of polynomials which are positive semidefinite on the set $S$. Obviously, the functions $f_i$ are all positive. Moreover, the product of two positive functions is also positive. Thus, taking all possible products of the functions $f_i$, we arrive at the Monoid - a set of functions which are positive on $S$.

\defnn{Multiplicative Monoid}{The multiplicative monoid $\mathcal{M}$ generated by elements $\{f_1,...,f_n\}\in\mathbf{R}[x]$ is the set
\[\mathcal{M}:=\left \{p\in\mathbf{R}[x]:p = \prod_{i=1}^n f_i^{k_i},\;k_i\in \mathbb{N}^0\right \}\]
Thus, $\mathcal{M}$ is the set of finite products of $\{f_1,...,f_n\}$.}

Now, we can add to the monoid a set of functions which are positive everywhere - namely the set of SOS polynomials. Taking all finite products of elements of these two sets, we get a larger set of polynomials which are positive on $S$ - The Cone (not to be confused with other traditional mathematical definitions of cone)

\defnn{Cone}{For given elements $\{f_1,...,f_n\}\in\mathbf{R}[x]$, the {\bf cone} $\mathcal{P}$ generated by $\{f_1,...,f_n\}\in\mathbf{R}[x]$ is the subset of $\mathbf{R}[x]$  defined as
\[\mathcal{P}:=\left\{p\in\mathbf{R}[x]: p=s_0 + \sum_{m \in \mathcal{M}} s_m m, \;,\;s_m \in \Sigma_x\right\}.\]
}

For a computationally simpler, yet equivalent definition of the {\bf cone}, let $\bar{\mathcal{M}}\subset\mathcal{M}$ be the set of products defined by
\[\bar{\mathcal{M}}:=\left\{p\in \mathbf{R}[x]: p = \prod_{i=1}^n f_i^{k_i}, k_i \in \{0,1\}\right\}\]
and let $M$ denote the cardinality of $\bar{\mathcal{M}}$.  Then the {\bf cone} $\mathcal{P}$ can be represented as
\[
\mathcal{P}:=\left\{p\in\mathbf{R}[x]: p=s_0 + \sum_{i=1}^M s_i m_i, \; m_i\in \bar{\mathcal{M}},\;s_i \in \Sigma_x\right\}.
\]

\noindent Note that $\mathcal{P}$ satisfies the following properties:
\begin{enumerate}
\item $a,b\in \mathcal{P}$ implies $a+b\in \mathcal{P}$
\item $a,b\in \mathcal{P}$ implies $a\cdot b\in \mathcal{P}$
\item $a\in \mathbf{R}[x_1,...,x_n]$ implies $a^2\in \mathcal{P}$
\end{enumerate}
\noindent\emph{Remark: }The {\bf cone} generated by $\{\emptyset\}$ is the cone of sum-of-squares polynomials, $\Sigma_x$.

\noindent The {\bf cone} $\mathcal{P}$ is not the largest set of polynomials which are positive semidefinite on $\mathcal{P}$. However, it is the largest set of such polynomials which is readily parameterized using positive matrices via SOS. Moreover, the Positivstellensatz tells us that the {\bf cone} is, in some sense, equivalent to the set of polynomial strictly positive on $S$. Before presenting this result, we consider the set of functions which are zero on $S$ - starting with the equality constraints $h_i(x)=0$. As with the monoid, all products of the functions $h_i$ are zero on $S$. Furthermore, the product of any function which is zero on $S$ with any other function is also zero on $S$. Thus the \emph{Ideal} is the set of products of the $h_i$ with arbitrary polynomials.

\defnn{Ideal}{The Ideal $\mathcal{I}$  generated by $\{h_1,...,h_n\}\in\mathbf{R}[x]$ is defined as
\[\mathcal{I} :=\left \{p\in\mathbf{R}[x]: p = \sum_{i=1}^n q_i h_i,\;q_i\in\mathbf{R}[x]\right \}\]
Note that $\mathcal{I}$ satisfies
\begin{enumerate}
\item $a,b\in \mathcal{I}$ implies $a+b\in \mathcal{I}$
\item $a\in \mathcal{I};b\in\mathbf{R}[x]$ implies $ ab\in \mathcal{I}$
\end{enumerate}
}

Intuitively, the ideal generated by a collection of polynomials is the set of polynomials that vanish when \emph{all} of the generating polynomials vanish. The following Positivstellensatz says that by combining the {\bf cone} and the {\bf ideal}, we get all polynomials which are positive on $S$.

\theoremnC{Stengle's Positivstellensatz}{
Given polynomials $\{f_1,f_2,...,f_{n_1}\}\subset \mathbf{R}[x]$, $\{g_1,g_2,...,g_{n_2}\}\subset \mathbf{R}[x]$, and $\{h_1,h_2,...,h_{n_3}\}\subset \mathbf{R}[x]$, let $\mathcal{P}$ be the cone generated by $\{f_i\}_{i=1,2,...,n_1}$, $\mathcal{M}$ be the multiplicative monoid generated by $\{g_j\}_{j=1,2,...,n_2}$, and $\mathcal{I}$ be the Ideal generated by $\{h_k\}_{k=1,2,...,n_3}$. Then, the following statements are equivalent:
\begin{enumerate}
\item $\{x\in\Realsn{n}\;:\;f_i(x)\geq 0, g_j(x) \neq 0, h_k(x) = 0, i=1,...n_1, j=1,...,n_2, k = 1,...,n_3\} = \emptyset$
\item  $\exists f\in \mathcal{P},\;\exists{} g\in\mathcal{M},\;\exists h\in\mathcal{I}$ s.t.
\[f + g^2 + h \equiv 0\]
\end{enumerate}
}
By noting that $p(x) >0$ for all $x \in S$ if and only if
\aligneq{
&\{x\in\Realsn{n}\;:\;-p(x)\geq 0,\; f_i(x)\geq 0, h_j(x) = 0, i=1,...n_1,\\ &\qquad\qquad\qquad\qquad\qquad\qquad\qquad\qquad j=1,...,n_3\} = \emptyset
\},
}
this theorem has the direct corollary
\theoremnC{Stengle's Corollary}{
Given polynomials $\{f_1,f_2,...,f_{n_1}\}\subset \mathbf{R}[x]$  and $\{h_1,h_2,...,h_{n_3}\}\subset \mathbf{R}[x]$, the following statements are equivalent:
\begin{enumerate}
\item $f(x)>0$ for all $x \in S:=\left \{x\in\Realsn{n}: f_i(x) \geq 0,\right.$\\$\left.i=1,...,n_1,\;h_i(x) = 0,i=1,...,n_3 \right \} $.
\item  There exist $s_i \in \Sigma_s$ and $q_i\in\mathbf{R}[x]$ such that
\aligneq{
p +  p\sum_{\substack{k \in \N^{n_1} \\ \norm{k}_\infty \le 1}} s_k \prod_{i=1}^{n_1} f_i^{k_i}  =& s_0  + \sum_{\substack{k \in \N^{n_1} \\ \norm{k}_\infty \le 1}} s_k \prod_{i=1}^{n_1} f_i^{k_i}\\
&+\sum_{i=1}^{n_3} q_i h_i }
\end{enumerate}
}

While Stengle's Positivstellensatz is interesting, the certificate of positivity it provides is bilinear in $p$ and $s_i$. If the set $S$ compact, Schmudgen's Positivstellensatz~\cite{Schmudgen} says we can neglect the summation on the left-hand side of the equation.
\theoremnC{Schmudgen's Corollary}{
Given polynomials $\{f_1,f_2,...,f_{n_1}\}\subset \mathbf{R}[x]$  and $\{h_1,h_2,...,h_{n_3}\}\subset \mathbf{R}[x]$, suppose 
\aligneq{
&S:=\left \{x\in\Realsn{n}: f_i(x) \geq 0,i=1,...,n_1,\;h_i(x) = 0,\right.\\
&\qquad\qquad\qquad\qquad\qquad\qquad\qquad\quad \left.i=1,...,n_3 \right \}
}
 is compact. Then the following are equivalent:
\begin{enumerate}
\item $f(x)>0$ for all $x \in S$.
\item  There exist $s_i \in \Sigma_s$ and $q_i\in\mathbf{R}[x]$ such that
\[p = s_0  + \sum_{\substack{k \in \N^{n_1} \\ \norm{k}_\infty \le 1}} s_k \prod_{i=1}^{n_1} f_i^{k_i} +\sum_{i=1}^{n_3} q_i h_i \]
\end{enumerate}
}
Now we have a linear parameterizations of the set of polynomials positive on $S$. However, the number of bases is rather large. If $S$ satisfies additional conditions~\cite{Putinar}, then most of the terms on the right-hand side of the equation can also be eliminated, leaving only
\[p = s_0  + \sum_{i=1}^{n_1} s_i f_i +\sum_{i=1}^{n_3} q_i h_i. \]
It is this representation we will use in this paper to parameterize the set of polynomials which are positive on a semialgebraic set.
Thus, by treating the $s_i$ and $q_i$ as variables, we can use convex Sum-of-Squares optimization algorithms (which convert the problem to an LMI) to search over the set of polynomials which are positive on a semialgebraic set in polynomial-time. Specifically, we use the Positivstellensatz to construct Lyapunov functions which are positive on bounded sets (see Sections 3 and 4).\par\noindent
For further details and proofs, we refer to \cite{Stengle} and \cite{ParriloThs}.

Note that the Positivstellensatz can also be thought of as a generalization of the S-procedure (as described in \cite{lmibook}). However, while the S-procedure certifies positivity of quadratic forms such that other quadratic forms are also positive, the Positivstellensatz can be used to obtain certificates of positivity for polynomials of arbitrary degree over semialgebraic sets.

\vspace{-.1in}
\subsection{Hybrid Systems}\label{subsection:hybsys}
\vspace{-.1in}
In this section, we present the formal definition of hybrid systems and executions that will be used in this paper. This framework is similar to the one used in, e.g.~\cite{HybVDSSchu}.

\defnn{Hybrid System}{
A hybrid system $H$ is a tuple:
\[H = (Q,E,D,F,G,R)\]
where
\begin{itemize}
\item $Q$ is a finite collection of discrete modes, states or indices.
\item $E\subset Q\times Q$ is a collection of edges.
\item $D = \{D_q\}_{q\in Q}$ is the collection of Domains associated with the discrete states, where for each $q\in Q$, $D_q\subseteq \Realsn{n}$.
\item $F = \{f_q\}_{q\in Q}$ is the collection of vector fields associated with the discrete states, where for each $q\in Q$, $f_q:D_q \rightarrow \Realsn{n}$.
\item $G = \{G_{e}\}_{e\in E}$ is a collection of guard sets, each associated with an edge. where for each $e = (q,q') \in E$, $G_{e} \subset D_{q}$
\item $R = \{\phi_{e}\}_{e\in E}$ is a collection of Reset Maps,  where for each $e=(q,q') E$, $\phi_{e}:G_{e}\rightarrow D_{q'}$.
\end{itemize}
}
Note that we also define the start and end functions $s,t:Q \times Q \rightarrow Q$ which act on the edges and indicate the start or end of that edge, so that for $e = (q,q')$, $s(e) = q$ and $t(e) = q'$.

\defnc{A cyclic hybrid system $H_c$ is a hybrid system where for each discrete mode $q \in Q$, there exists a unique edge $e_q\in E$ such that $s(e_q)=q$ and a unique edge $e_q' \in E$ such that $q = t(e_q')$ and such that set of edges and modes forms a connected digraph.}
In a cyclic hybrid system, each discrete mode is the source of only one edge, and the target of only one edge. Let $e(q)=e_q$, then $q=t(e(s(\cdots e(t(e(q))))))$ - the sequence of edges will eventually return to the original mode.

\assn{
\noindent In this paper, we consider hybrid systems with polynomial vector fields and resets, and semialgebraic domains and guard sets. This implies that for every hybrid system, we there exists a set of polynomials $g_{q,i}$, $h_{e,k}$ for $q \in Q$, $e \in E$, $i=1,\cdots,K_{q}$ and $k = 1,\cdots,N_{q}$ for some $K_q,N_q>0$ such that
\begin{equation}
D_q = \{x\in \Realsn{n}\;: \;g_{q,i}(x)\geq 0,\; i=1,2,\cdots,K_{q}\}
\end{equation}
and
\begin{eqnarray}
\nonumber G_e = \{x\in \Realsn{n}\;:\; h_{e,0}(x) = 0,\; h_{e,k}(x)\geq 0,\qquad\\
\;\qquad\qquad\qquad\qquad\qquad k = 1,2,\cdots,N_{q} \}
\end{eqnarray}
Furthermore, this implies that for each $e=(q,q')\in E$, there exist polynomials $\phi_{e,j} \in \mathbf{R}[x]$ such that the reset map $\phi_{e}$ is given by the vector-valued polynomial function
\begin{equation}
\phi_{e} = [\phi_{e,1},\cdots,\phi_{e,n}]^T.
\end{equation}
}

The Cauchy problem of defining solutions for hybrid systems is defined in terms of an execution.

\defnn{Hybrid System Execution}{We say that the tuple
\[\chi = (I,T,p,C)\]
where
\begin{itemize}
\item $I\subseteq \mathbb{N}$ indexes the intervals of time on which the trajectory continuously evolves.
\item $T = \{T_i\}_{i\in I}$ are a set of open time intervals associated with points in time $\tau_i$ as $T_i = (\tau_i,\tau_{i+1}) \subset \Realsn{n}_+$ where $T_{i+1} = (\tau_{i+1},\tau_{i+2})$.
\item $p:I\rightarrow Q$ maps each interval to a discrete mode.
\item $C = \{c_{i}\}_{i\in I}$ is a set of continuously differentiable functions where $c_i \in \mathcal{C}[T_i]$.
\end{itemize}
is an \emph{execution} of the hybrid system $H=F(Q,E,D,F,G,R)$ with initial condition $(q_0,x_0)$ if
\begin{enumerate}
\item $c_{1}(0) = x_0$ and $p(1) = q_0$.
\item $\dot{c}_{i}(t)= f_{p(i)}(c_{i}(t))$ for $t\in T_i$ for every $i \in I$.
\item $c_i(t)\in D_{p(i)}$ for $t\in T_i$ for every $i \in I$.
\item $c_i(\tau_{i+1})\in G_{(p(i),p(i+1))}$ for every $i\in I$.
\item $c_{i+1}(\tau_{i+1}) = \phi_{(p(i),p(i+1))}(c_i(\tau_i))$ for every $i\in I$.
\end{enumerate}
}

Note that an execution does not require $\lim_{i\rightarrow \infty} \tau_i = \infty$, so the solution may not be defined for all time.
\vspace{-.1in}
\subsection{Zeno Stability in Hybrid Dynamical Systems}\label{subsection:zeno}
\vspace{-.1in}
In this section, we define Zeno executions, Zeno equilibria, and Zeno stability. We also present a Lyapunov theorem for Zeno stability given in~\cite{AmesLyap} and~\cite{TeelGoebLyap}.
\defnn{Zeno Execution}{
We say an execution $\chi = (I,T,p,C)$ starting from $(q_0,x_0)$ of a hybrid System $=(Q,E,D,F,G,R)$ is Zeno if
\begin{enumerate}
\item $I = \mathbb{N}$
\item $\lim_{i \rightarrow \infty} \tau_i < \infty$
\end{enumerate}
}
Thus a Zeno execution is an execution which undergoes infinite discrete transitions in finite-time.
\defnn{Zeno Equilibrium}{A set $z = \{z_q\}_{q\in Q}$ with $z_q \in D_q$ is a Zeno equilibrium of a Hybrid System $H=(Q,E,D,F,G,R)$ if it satisfies
\begin{enumerate}
\item For each edge $e = (q,q') \in E$, $z_q\in G_e$ and $\phi_e(z_q) = z_{q'}$.
\item $f_q(z_q)\neq 0$ for all $q \in Q$.
\end{enumerate}
For any $z\in \{z_q\}_{q\in Q}$, where $\{z_q\}_{q\in Q}$ is a Zeno equilibrium of a cyclic hybrid system $H_c$,
\[\left (\phi_{i-1}\circ \cdots \circ \phi_0 \cdots \phi_i\right ) (z) = z\]
}
By definition, a Zeno equilibrium is NOT an equilibrium point ($f_q(z_q)\neq0$). Although the results of this paper may be readily extended to consider classical (non-Zeno) stability, such results already exist in the literature. Note that a Zeno equilibrium also defines a Zeno execution with $\lim_{i \rightarrow \infty} \tau_i=0$.
A Zeno equilibrium is \emph{isolated} if there exists neighborhoods $X_q$ of $z_q$ such that for any other Zeno equilibrium $\hat z$, $\hat z_q \not \in X_q$ for some $q \in Q$. That is, the equilibrium is strictly separated from other equilibria.

\defnn{Zeno Stability}{%
Let \[H = (Q,E,D,F,G,R)\] be a hybrid system, and let $z = \{z_q\}_{q\in Q}$ be a Zeno equilibrium.  The set $z$ is Zeno stable if, for each $q\in Q$, there exist neighborhoods $Z_q$, where $z_q\in Z_q$, such that for any initial condition $(x_0,q_0)\in \bigcup_{q\in Q}(Z_q,q)$, the execution $\chi = (I,T,p,C)$, with initial condition $(x_0,q_0)$ is Zeno, and for any $\epsilon$, there exists an $N \in \N$ such that $i>N$ implies $\norm{c_i(T_i(2))- z_{p(i)}}\le \epsilon$.
}
We give a slight variation of the conditions for Zeno stability of cyclic hybrid systems presented in~\cite{AmesLyap,TeelGoebLyap}.\par

\theoremnC{Lamperski and Ames}{\label{thm2_v1}
\noindent Consider a cyclic hybrid system $H = (Q,E,D,F,G,R)$, with an isolated Zeno equilibrium $\{z_q\}_{q\in Q}$. Let $\{W_q\}_{q\in Q}$ be a collection of open neighborhoods of $\{z_q\}_{q\in Q}$. Suppose there exist continuously differentiable functions $V_q:\Realsn{n}\rightarrow \Reals$ and $B_q:\Realsn{n}\rightarrow \Reals$, and constants $r_q\in [0,1]$, $\gamma_a,\gamma_b \ge 0$, for every $q \in Q$ where $r_q<1$ for some $q$ and such that
\begin{align}
V_q(x) &> 0 \quad \text{for all } x\in W_q\backslash z_q, q \in Q\label{eq:EC1}\\
V_q(z_q) &= 0, \quad \text{for all } q \in Q\;\label{eq:EC1b}\\
\grad {V}_q^T(x)f_q(x) &\leq 0\quad \text{for all } x\in W_q, \, q \in Q\label{eq:EC2}\\
B_q(x) &\geq 0\quad \text{for all } x\in W_q, \, q \in Q \label{eq:EC3}\\
\grad {B}_q^T(x)f_q(x) &< 0\quad \text{for all } x\in W_q, \, q \in Q  \label{eq:EC4}\\
\nonumber V_{q'}(\phi_{e}(x)) &\leq r_q V_q(x),\;\;\text{for all } x\in G_{e}\cap W_q  \\&\qquad\qquad\text{ and all }e = (q,q') \in E \label{eq:C1}\\
\nonumber B_q(\phi_{e}(x)) &\leq \gamma_b \left(V_q(\phi_{e}(x)) \right)^{\gamma_a}, \\\text{for all }x \in G_{e} &\cap W_q \text{ and all }e =(q,q') \in E.\label{eq:C2}\vspace{-8mm}
\end{align}
Then $\{z_q\}_{q\in Q}$ is Zeno stable.
}
\par
\noindent \par
 In this paper, we use a simplified version of \ref{thm2_v1}. Although we have implemented and tested the conditions in Theorem~\ref{thm2_v1}, numerical tests indicate little or no additional conservativity is implied by using the following simplification. \par
\theoremC{\label{thm2}
Let $H = (Q,E,D,F,G,R)$ be a cyclic hybrid system, and let $z = \{z_q\}_{q\in Q}$ be a Zeno equilibrium. Suppose we have $\{W_q\}_{q\in Q}$ with $W_q\subset D_q$ and $W_q$ a neighborhood of $z_q$ for each $q \in Q$. Now suppose that there exist continuously differentiable functions $V_q:W_q\rightarrow \Reals$, and constants $r_q\in (0,1]$, $\gamma>0$ for $q \in Q$ where $r_q<1$ for some $q\in Q$ and such that
\begin{flalign}
V_q(x) &> 0 \quad \text{for all } x\in W_q\backslash z_q, q \in Q \label{eq:NC1}\\
V_q(z_q) &= 0, \quad \text{for all } q \in Q \label{eq:NC2}\\
\grad {V}_q^T(x)f_q(x) &\leq-\gamma\;\quad \text{for all } x\in W_q, \, q \in Q \label{eq:NC3}\\
 \nonumber V_{q'}(\phi_{e}(x))&\leq r_q V_q(x) ,\;\text{for all } x \in G_{e} \cap W_q\\&\qquad\quad\text{ and all }e =(q,q') \in E . \label{eq:NC4}
\end{flalign}
then $z$ is Zeno stable.
}
\vspace{-.4in}
\ProofC{
We show that if for each $q\in Q$, we can find a $V_q$ such that (\ref{eq:NC1})-(\ref{eq:NC4}) are satisfied, then the same $V_q$ also satisfies (\ref{eq:EC1})-(\ref{eq:C2}). From inspection, it is clear that if $V_q$ satisfies (\ref{eq:NC1})-(\ref{eq:NC4}), then (\ref{eq:EC1})-(\ref{eq:EC2}) and (\ref{eq:C1}) are satisfied. Second, choose $B_q=V_q$ for each $q\in Q$. From inspection, it is clear that $B_q$ also satisfies (\ref{eq:EC3}) and (\ref{eq:EC4}). Last, if $\gamma_a = \gamma_b=1$, we get $V_q\leq V_q$, where the equality holds. From this, we see that for each $q\in Q$, $V_q$ also satisfies (\ref{eq:C2}).
}
\vspace{-.95in}
\section{Using Sum-of-Squares Programming to prove Zeno Stability}
\vspace{-.1in}
Theorem~\ref{thm2} provides sufficient conditions for Zeno stability in cyclic hybrid systems. We now show how to enforce these conditions using sum-of-squares programming.\par
Let $H = (Q,E,D,F,G,R)$ be a hybrid system, and let $z=\{z_q\}_{q\in Q}$. Let $\{W_q\}_{q\in Q}$ be a collection of neighborhoods of $\{z_q\}_{q\in Q}$. Suppose that each $W_q$ is a semialgebraic set defined as
\[
W_q:=\{x\in D_q\, : \, w_{qk}(x)\ge 0, k=1,2,...,K_{qw}\}
\]
where $w_{qk}\in\mathbf{R}[x]$. For example, if $w_{q1}(x) =1-x^Tx$, then $W_q$ is the unit ball intersected with $D_q$. \\ \noindent
We define feasibility problem 1:\\\\
\noindent{\bf Feasibility Problem 1:}\\
For hybrid system $H = (Q,E,D,F,G,R)$, find
\begin{itemize}
\item $a_{qk}$, $c_{qk}$, $i_{qk}$, $\in \Sigma_x$, for $k=1,2,...,K_{qw}$ and $q \in Q$;
\item $b_{qk}$, $d_{qk}$, $j_{qk} \in \Sigma_x$, for $k=1,2,...,K_q$ and $q \in Q$.
\item $m_{e,l} \in \Sigma_x$ for $e \in E$ and $l=1,2,...,N_q$
\item $V_q$, $m_{e,0}\in \mathbf{R}[x]$ for $e \in E$ and $q \in Q$.
\item Constants $\alpha, \gamma > 0$, and $r_q\in (0,1]$ for $q \in Q$ such that $r_q<1$ for some $q \in Q$.
\end{itemize}
such that
\begin{flalign}
\nonumber &V_q - \alpha x^T x-\sum_{k=1}^{K_{qw}} a_{qk}w_{qk} -\sum_{k=1}^{K_q} b_{qk}g_{qk} \in \Sigma_x \quad \\\label{eq:mVp1}&\qquad\qquad\qquad\qquad\qquad\qquad\qquad\text{for all }q\in Q\\
\label{eq:mVz}& V_q(z_q) = 0 \quad \text{ for all }q \in Q \\
\nonumber&-\grad V_q^Tf_q - \gamma -\sum_{k=1}^{K_{qw}} c_{qk}w_{qk} -\sum_{k=1}^{K_q} d_{qk}g_{qk} \in \Sigma_x\\\label{eq:mdVn1}&\qquad\qquad\qquad\qquad\qquad\qquad \qquad \text{ for all } q \in Q\\
&\nonumber r_q V_q - V_{q'}(\phi_{e})-  m_{e,0}h_{e,0} - \sum_{l=1}^{N_q} m_{e,l}h_{e,l}\\&\qquad\qquad\quad\;\;\;- \sum_{k=1}^{K_{qw}} i_{qk}w_{qk} - \sum_{k=1}^{K_{q}} j_{qk}g_{qk} \in \Sigma_{x} \nonumber\\&\qquad\qquad\qquad\qquad\qquad\text{ for all } e = (q,q') \in E \label{eq:mVC1}
\end{flalign}
\noindent\theoremC{\label{thm_SOS}
Consider a hybrid system \[H = (Q,E,D,F,G,R),\]
 and let $z=\{z_q\}_{q\in Q}$. If Feasibility Problem 1 has a solution, then $z$ is Zeno stable.
}
\ProofC{
To prove the theorem we show that if $V_q$, $q \in Q$ are elements of a solution of Feasibility Problem 1, then for each $q\in Q$, the same $V_q$ also satisfy (\ref{eq:EC1})-(\ref{eq:C1}) of Theorem~\ref{thm2_v1}.
That is, we show that if the $V_q$ satisfy (\ref{eq:mVp1})-(\ref{eq:mVC1}), then the same $V_q$ also satisfies (\ref{eq:NC1})-(\ref{eq:NC4}).  \\
First, we observe that~\eqref{eq:mVz} directly implies~\eqref{eq:NC2}. Next, from~\eqref{eq:mVp1}, we know that
\begin{align*}
& V_q(x) \geq \sum_{k=1}^{K_{qw}} a_{qk}(x)w_{qk}(x) +\sum_{k=1}^{K_q}b_{qk}(x)g_{qk}(x) + \alpha x^T x.
\end{align*}
Since $a_{qk}(x)$ and $b_{qk}(x)$ are SOS they are nonnegative. Furthermore, by the definitions of $W_q$ and $D_q$, we know $w_{qk}(x)$ and $g_{qk}(x)$ are non-negative on $W_q$. Thus $V_q(x) \geq \alpha x^T x$ for all $x\in W_q\subset D_q$. Thus, (\ref{eq:mVp1}) implies (\ref{eq:NC1}) is satisfied. Similarly, from (\ref{eq:mdVn1}),
\aligneq{
 -\grad V_q^T(x)f_q(x) -\gamma\geq& \sum_{k=1}^{K_{qw}} c_{qk}(x)w_{qk}(x)\\
&+\sum_{k=1}^{K_q} d_{qk}(x)g_{qk}(x).
}
As before, $c_{qk}(x)$ and $d_{qk}(x)$ are SOS and hence $\grad V_q(x) ^T f_q(x) \le -\gamma$ for $x\in W_q$ which implies (\ref{eq:NC3}) is satisfied. Next, from~\eqref{eq:mVC1} we have that for all $e = (q,q') \in Q$,
\begin{align*}
& r_qV_q(x) - V_{q'}(\phi_{e}(x))\\ &\geq m_{e,0}(x)h_{e,0}(x) 
 +\hspace{-.14cm} \sum_{l=1}^{N_q} m_{e,l}(x)h_{e,l}(x)+\hspace{-.14cm}\sum_{k=1}^{K_{qw}} i_{qk}(x)w_{qk}(x) \\& \qquad\qquad\qquad\qquad\qquad\qquad\qquad\qquad+ \hspace{-.14cm}\sum_{k=1}^{K_q} j_{qk}(x)g_{qk}(x).
\end{align*}
First note that $h_{e,0}(x) = 0$ and hence $m_{e,0}(x)h_{e,0}(x) = 0$ on $G_e$. Since $m_{e,l} \in \Sigma_x$, we have  $m_{e,l}(x)h_{e,l}(x)\ge 0$ on $G_e$. As before, $j_{qk}(x)g_{qk}(x)$ and $i_{qk}(x)w_{qk}(x)$ are non-negative on $W_q$. It follows that $r_q V_q(x) - V_{q'}(\phi_{e}(x)) \geq 0$ when $x\in G_e \cap W_q $ for all $e = (q,q') \in E$. Thus (\ref{eq:mVC1}) implies (\ref{eq:NC4}). \par
\noindent  Thus we conclude that any solution $\{V_q\}_{q \in Q}$ of Feasibility Problem 1 satisfies the conditions~(\ref{eq:NC1})-(\ref{eq:NC4}) of Theorem~\ref{thm2} which by Theorem~\ref{thm2_v1} implies Zeno stability of $z$.
}%

\section{Examples}
\begin{figure}[h]
\centering
\includegraphics[width=0.8\columnwidth]{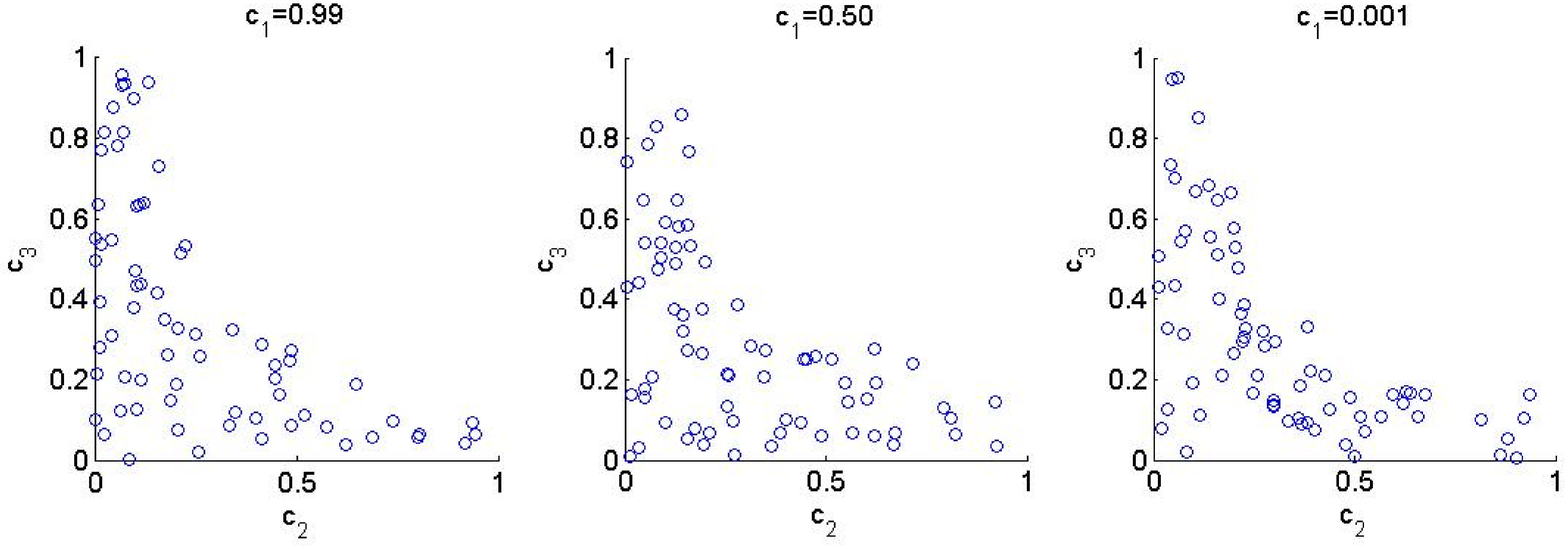}
\caption{Values of $c_2$ and $c_3$ for fixed $c_1$}
\label{fig:3}
\end{figure}

\begin{figure}[h]
\centering
\includegraphics[width=0.8\columnwidth]{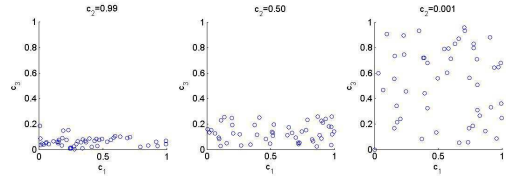}
\caption{Values of $c_1$ and $c_3$ for fixed $c_2$}
\label{fig:4}
\end{figure}

\begin{figure}[h]
\centering
\includegraphics[width=0.8\columnwidth]{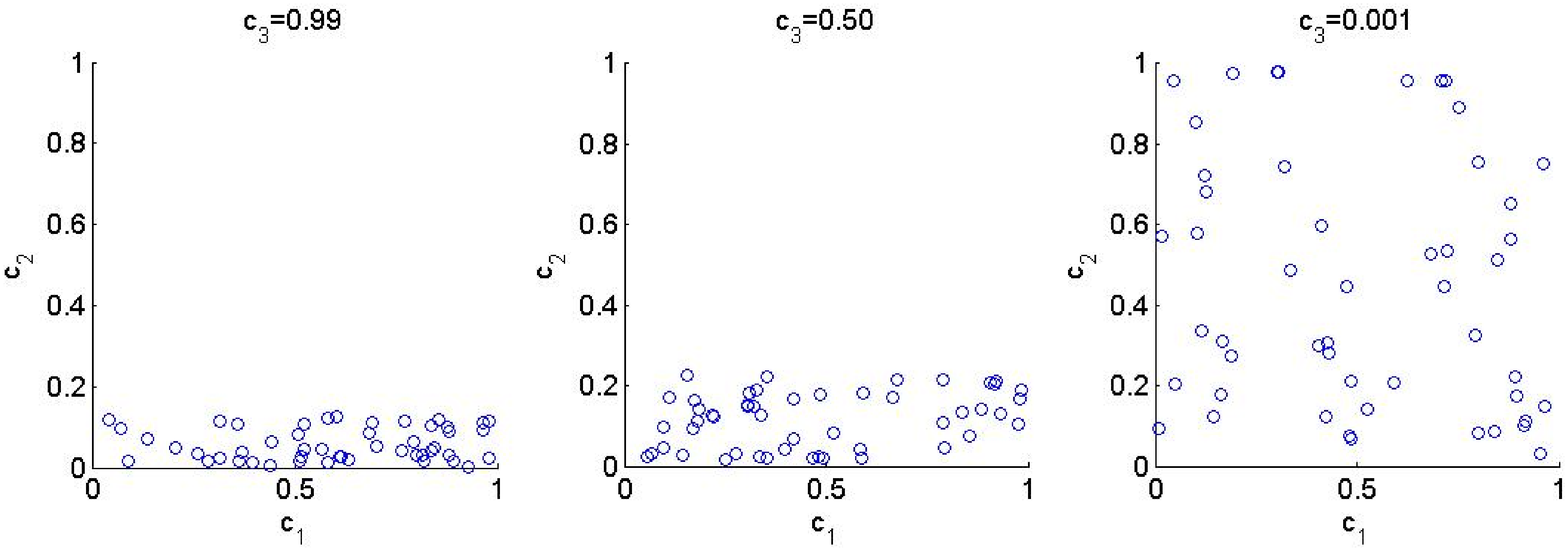}
\caption{Values of $c_1$ and $c_2$  for fixed $c_3$}
\label{fig:5}
\end{figure}

\vspace{-.1in}
In this section, we show how the proposed method can be applied to some simple examples.
\exmpl{\textbf{(Bouncing Ball)} Define the nonlinear hybrid system $H_N$ as:
\[H_N = (Q,E,D,F,G,R)\]
where
\begin{itemize}
\item $Q = \{1\}$
\item $E = \{(1,1)\}$
\item $D := \{x\in\Realsn{2}: x_1\geq 0\}$
\item $G := \{x\in\Realsn{2}: x_1=0,\;x_2\leq 0\}$
\item $F = \{f\}$, where
\[\dot{x} = f(x) = \rmatrix{c}{x_2\\-g+c_1x_2^2}\]
\item $R= \phi(x) = [0,-c_2x_2(1-c_3x_2^2)]^T$. Here, $c_1$, $c_2$, and $c_3$ can be any positive constants satisfying $c_i < 1$.
\end{itemize}
}
{\bf Results}
\begin{figure}[h]
\centering
\includegraphics[width=0.75\columnwidth]{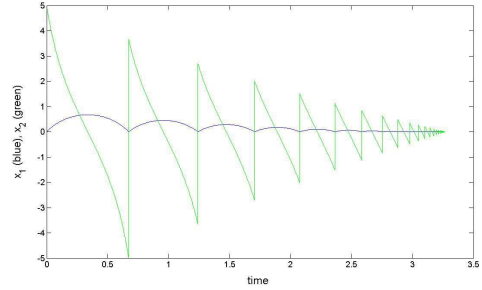}
\caption{Nonlinear Hybrid System with $c_1=0.5$, $c_2=0.8$, $c_3=0.001$}
\label{fig:1}
\end{figure}
Our goal is to analyze Zeno stability of the $z = (0,0)$ Zeno equilibrium. We used SOSTOOLS to search for a 6{th}-order (degree 6 polynomial) $V(x)$ and associated SOS multipliers satisfying the conditions of Feasibility Problem 1. The neighborhoods we chose were $W_q:=\{x\in D_q \, : \, w(x) = 25-x^T x \geq 0\}$ - which is the ball of radius 5. We were able to show Zeno stability for a range of parameters $c_i$. A numerical simulation of the system is shown in Figure~\ref{fig:1} for $c_1=0.5$, $c_2=0.8$, $c_3=0.001$. To better illustrate the range of Zeno-stable parameters, we used a Monte-Carlo approach to selection of the parameters $c_i$. At each set of parameters, the algorithm was able to prove stability or return a certificate of infeasibility. The results are seen in Figures~\ref{fig:3} -~\ref{fig:5}. In Figure~\ref{fig:3}, we estimate the set of Zeno-stable values of $c_2$ and $c_3$ for three different values of $c_1$. In Figure~\ref{fig:4}, we estimate the set of Zeno-stable values of $c_1$ and $c_3$ for three different values of $c_2$. Finally, in Figure~\ref{fig:5}, we estimate the set of Zeno-stable values of $c_1$ and $c_2$ for three different values of $c_3$.
We note from Figure~\ref{fig:3} that the range of values of $c_1$ and $c_2$ for which $z$ is Zeno stable does not seem to depend on $c_1$.
\exmpl{
\textbf{(Sliding Mode Control)} We consider the hybrid system $H = (Q,E,D,F,G,R)$ where
\begin{itemize}
\item $Q = \{1,2\}$
\item $E = \{(1,2),(2,1)\}$
\item $D = \{D_1,D_2\}$ where
\aligneq{
D_1 := \{x\in \Realsn{2}: x_1+x_2\geq 0\}\\
D_2 := \{x\in \Realsn{2}: x_1+x_2\leq 0\}
}
\item $F = \{f_1,f_2\}$ where
\aligneq{
f_1 = \rmatrix{c}{x_2\\
3(x_2^2+x_1^2)}\\
f_2 = \rmatrix{c}{x_2\\
-(x_2^2+x_1^2)}}
\item $G = \{G_{12},G_{21}\}$ where
\aligneq{
G_{12}=G_{21} := \{x\in \Realsn{2}: x_1+x_2 = 0\}
}
\item $R= \{\phi_{12}(x),\phi_{21}(x)\}$ where each $\phi_{ij}(x) = x$.
\end{itemize}

\begin{figure}[h]
\centering
\includegraphics[width=0.6\columnwidth]{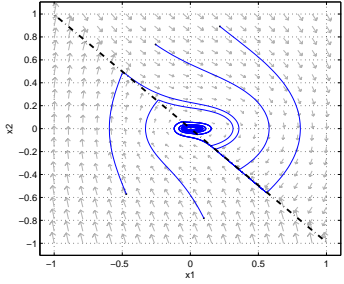}
\caption{Hybrid system of Example 2. The dashed line indicates $x_2+x_1=0$.}
\label{fig:slidingmode}
\end{figure}

{\bf Results:}
\noindent For a slightly modified form of sliding-mode controller, we examined stability of the Zeno equilibrium
\[ z = \left\{\bmat{0\\0},\bmat{0\\0}\right\}.\]
Note that this is actually a true equilibrium. However, as mentioned, these tools also apply to such classical problems. For our analysis, we analyze Zeno stability in the neighborhoods $W_{1}:= \{x\in D_1: |x|\le 1\}$ and $W_{2}:= \{x\in D_2: |x|\le 1\}$. We used SOSTOOLS to find that  verification of the Conditions of Feasibility Problem 1 required the use of degree 8 polynomials. Naturally, the polynomials are too long for publication. However, a simulation illustrating stability is shown in Figure~\ref{fig:slidingmode}.
}

\exmpl{
\textbf{(Gain-Scheduling)} Consider the hybrid system $H = (Q,E,D,F,G,R)$, where
\begin{itemize}
\item $Q = \{1,2,3\}$
\item $E = \{(1,2),(2,3),(3,1)\}$
\item $D := \{D_1,D_2,D_3\}$ where
\aligneq{
&D_1 = \{x\in\Realsn{2}: x_1>0, x_2 + \frac{1}{2}x_1 \geq  0 \}\\
&D_2 = \{x\in\Realsn{2}: x_2-\frac{1}{2}x_1\geq 0, x_2 + \frac{1}{2}x_1 <  0 \}\\
&D_3 = \{x\in\Realsn{2}: x_1<0, x_2 + \frac{1}{2}x_1 \geq  0 \}}
\item $F = \{f_1,f_2,f_3\}$, where
\aligneq{
\dot{x} &= f_1(x) = \rmatrix{c}{x_2\\-5x_1-x_2}\\
\dot{x} &= f_2(x) = \rmatrix{c}{-x_1^2-3\\2x_2^2-\frac{1}{2}x_1^2}\\
\dot{x} &= f_3(x) = \rmatrix{c}{x_2^2+x_1\\-3x_1}
}
\item $G := \{G_{12},G_{23}, G_{31}\}$ where
\aligneq{
&G_{12}:=\left\{x\in \Realsn{2}: x_2\leq 0, \frac{1}{2}x_1+x_2=0 \right\}\\
&G_{23}:=\left\{x\in \Realsn{2}: x_2\leq 0, \frac{1}{2}x_1-x_2=0 \right\}\\
&G_{31}:=\left\{x\in \Realsn{2}: x_2> 0, x_1=0 \right\}
}

\item $R= \{\phi_{12}(x),\phi_{23}(x),\phi_{31}(x)\}$ where each $\phi_{ij}(x) = x$.
\end{itemize}

\begin{figure}[h]
\centering
\includegraphics[width=0.6\columnwidth]{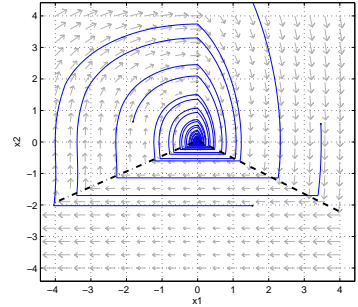}
\caption{Hybrid System in Example 3. Dashed line indicates $G_{12}$, dash-dotted line indicates $G_{23}$ and dotted line indicates $G_{31}$ }
\label{fig:nonlinear2}

\end{figure}

\noindent {\bf Results:} Zeno behavior such as exhibited by this system can arise due to, e.g. gain scheduling and may result in the state getting ``stuck'' at a non-equilibrium position. A phase portrait of the system is given in Figure~\ref{fig:nonlinear2}.

In this case, the equilibrium is Zeno and occurs at
\[ z = \left\{\bmat{0\\0},\,\bmat{0\\0},\,\bmat{0\\0}\right\}.
\]
We use the neighborhoods $W_{1}:= \{x\in D_1: |x|< 1\}$, $W_{3}:= \{x\in D_3: |x| \leq 1\}$, and  $W_{3}:= \{x\in D_3: |x| \leq 1\}$.
We were able to solve Feasibility Problem 1 for this system using degree 8 polynomials, implying Zeno stability.
}
\vspace{-.1in}
\section{Zeno Stability in Systems with Uncertainties}
\vspace{-.1in}
In this section, we show how the Sum-of-Squares Methodology can be leveraged to verify Zeno stability in cyclic hybrid systems with parametric uncertainty in the guard set, vector fields, and reset maps. To do this, we suppose the set of admissible uncertain parameters is a semialgebraic set of the form
\begin{equation}\label{eq:params}
\mathcal{P}:=\{p\in\Realsn{n_p}: \pi_{k}(p)\geq 0, k=1,2,...,K_1 \},
\end{equation}
where the $\pi_k \in \R[x]$.

In this paper, we use the following model of an uncertain hybrid system:\par

\noindent \textbf{Uncertain Hybrid Model:}
\noindent We consider a parameterized hybrid system $H(p) = (Q,E,D(p),F(p),G(p),R(p))$ where the set of domains $D(p) = \{D_q(p)\}_q$ is defined as
\begin{flalign}
\nonumber&D_q(p) := \{x\in \Realsn{n}\;: \;g_{qk}(x,p)\geq 0,\;\qquad\qquad\qquad\quad\\&\qquad\qquad\qquad\qquad\qquad k=1,2,\cdots,K_{q \in Q}\},
\end{flalign}
with $g_{qk}\in \mathbf{R}[x,p]$. The set of guard sets $G(p) = \{G_e(p)\}_{e \in E}$ are defined as
\begin{flalign}
\nonumber&G_e = \{x\in D_q(p)\;:\; h_{e,0}(x,p) = 0,\;\\ &\qquad\qquad\qquad h_{e,k}(x,p)\geq 0, k = 1,2,\cdots,N_{q} \},
\end{flalign}
with $h_{ek}\in\mathbf{R}[x,p]$. The set of reset maps $R(p):=\{\phi(p)\}_{e \in E}$ are defined by polynomials
\begin{equation}
\phi_{e}(x,p) = [\phi_{e,1}(x,p),\cdots,\phi_{e,n}(x,p)]^T
\end{equation}
where $\phi_{e,j} \in \mathbf{R}[x,p]$. The set of vector fields $F(p)=\{f_q(p)\}_{q \in Q}$ is likewise assumed to be a vector of polynomials.

Now, we present a parameterized version of Theorem~\ref{thm2} for this class of uncertain hybrid systems:
\theoremC{\label{thm2r}
Let $H(p) = (Q,E,D(p),F(p),G(p),R(p))$ be a set of cyclic hybrid systems parameterized by $p \in P:=\{p\in\Reals: \pi_{k}(p)\geq 0, k=1,2,...,K_1 \}$ and with common Zeno equilibrium $z = \{z_q\}_{q\in Q}$. Let $\{W_q\}_{q \in Q}$ be a collection of open neighborhoods of $\{z_q\}_{q\in Q}$.
Suppose that there exist continuously differentiable functions $V_q:D_q\times P \rightarrow \Reals$, and constants $r_q\in (0,1]$ for $q \in Q$ and $\gamma>0$, where $r_q<1$ for some $q\in Q$ and such that
\begin{flalign}
\label{eq:NC1r} &V_q(x,p) > 0 \quad \text{for all } x\in W_q\backslash z_q,\;p\in P,\; q \in Q ,\\
&V_q(z_q,p) = 0, \quad \text{for all } q \in Q,\; p\in P, \label{eq:NC2r}\\
\nonumber &\left ( \nabla_x {V}_q(x,p) \right )^T f_q(x)\leq-\gamma\;\quad \\&\qquad\qquad\qquad\quad\text{for all } x\in W_q, \, q \in Q ,\label{eq:NC3r}\\
 \nonumber &r_q V_q(x,p) \geq V_{q'}(\phi_{e}(x,p),p)\;  \\\label{eq:NC4r}&\qquad\text{for all } e =(q,q') \in E \text{ and }x \in G_{e} \cap W_q. 
\end{flalign}
Then $z$ is a Zeno stable Zeno equilibrium of $H(p)$  for any $p\in P$.
}

As before, consider neighborhoods $W_q$ of the form
\[W_q:=\{x\in\Realsn{n}: w_{qk}(x)>0, k=1,2,...,K_q\}\]
where $w_{qk}(x) \in \mathbf{R}[x]$. We now define a new SOS feasibility problem.\\\\
\noindent{\bf Feasibility Problem 2:\\}
For set of hybrid systems \vspace{-.15in}\[H(p) = (Q,E,D(p),F(p),G(p),R(p))\] defined above, find
\vspace{-.15in}
\begin{itemize}
\item $a_{qk}$, $c_{qk}$, $i_{qk}$, $\in \Sigma_{x,p}$, for $k=1,2,...,K_{qw}$, $p\in P$ and $q \in Q$;
\item $b_{qk}$, $d_{qk}$, $j_{qk} \in \Sigma_{x,p}$, for $k=1,2,...,K_q$, $p\in  P$, and $q \in Q$.
\item $ \eta_{qk}$, $ \beta_{qk}$, $ \zeta_{qk} \in \Sigma_{x,p}$, for $k=1,2,...,K_1$, $p\in P$, and $q \in Q$.
\item $m_{e,l} \in \Sigma_{x,p}$ for $e \in E$, $p\in P$, and $l=1,2,...,N_q$
\item $V_q$, $m_{e,0}\in \mathbf{R}[x,p]$ for $e \in E$, $p\in P$, and $q \in Q$.
\item Constants $\alpha, \gamma > 0$, $\{r_q\}_{q\in Q}\in (0,1]$ such that $r_q<1$ for some $q \in Q$.
\end{itemize}
such that
\begin{flalign}
& \nonumber V_q - \alpha x^T x-\sum_{k=1}^{K_{qw}} a_{qk}w_{qk} -\sum_{k=1}^{K_q} b_{qk}g_{qk}\\&\qquad\qquad\;\;\; - \sum_{k=1}^{K_1} \eta_{qk}\pi_{qk} \in \Sigma_{x,p} \quad \text{for all }q\in Q\label{eq:R1}\\
 &V_q(z_q,p) = 0 \quad \text{ for all }q \in Q \\
\nonumber &-\grad V_q^Tf_q - \gamma -\sum_{k=1}^{K_{qw}} c_{qk}w_{qk} -\sum_{k=1}^{K_q} d_{qk}g_{qk} \\&\qquad\qquad\;\;- \sum_{k=1}^{K_1} \beta_{qk}\pi_{qk} \in \Sigma_{x,p} \quad \text{ for all } q \in Q\label{eq:R2}\\
\nonumber &r_q V_q \hspace{-.25mm}-\hspace{-.25mm} V_{q'}(\phi_{e})\hspace{-.25mm}-\hspace{-.25mm}  m_{e,0}h_{e,0}\hspace{-.25mm} - \hspace{-.25mm}\sum_{l=1}^{N_q} m_{e,l}h_{e,l}\hspace{-.25mm}  \\\nonumber&\quad - \hspace{-.25mm}\sum_{k=1}^{K_{qw}} i_{qk}w_{qk} \hspace{-.25mm}-\sum_{k=1}^{K_{q}} j_{qk}g_{qk} \hspace{-.25mm}- \hspace{-.25mm}\sum_{k=1}^{K_1} \zeta_{qk}\pi_{qk} \in \Sigma_{x,p} \\&\qquad\qquad\qquad\qquad\qquad\quad \text{ for all } e = (q,q') \in E. \label{eq:R3}
\end{flalign}

\par\noindent
\theoremC{\label{robust}
Let $H(p) = (Q,E,D(p),F(p),G(p),R(p))$ be a set of cyclic hybrid systems parameterized by $p \in P:=\{p\in\Reals: \pi_{k}(p)\geq 0, k=1,2,...,K_1 \}$ and with common Zeno equilibrium $z = \{z_q\}_{q\in Q}$. If Feasibility Problem 2 has a solution, then $z$ is a Zeno stable Zeno equilibrium of $H(p)$  for any $p\in P$.
}
\ProofC{
The proof is similar to that of Theorem~\ref{thm_SOS}, except that is the case of parametric uncertainty, we have that $\pi_{k} (p) \geq 0$ for all $p\in P$ and implicitly all $x \in W_q$. This implies that $\eta_{qk}(x,p)\pi_{qk}(p)$, $\beta_{qk}(x,p)\pi_{qk}(p)$, and $\zeta_{qk}(x,p)\pi_{qk}(p)$ are also non-negative. Hence, by similar logic to that employed in the proof of Theorem~\ref{thm_SOS}, we have that the functions $V_q$ satisfy the Conditions (\ref{eq:NC1})-(\ref{eq:NC4}). By Theorem~\ref{thm2r}, this implies that $z$ is a Zeno stable Zeno equilibrium of $H(p)$  for any $p\in P$.}
\vspace{-.1in}
\subsection{Numerical Examples}
\vspace{-.1in}
We now present two examples which illustrate Theorem~\ref{robust}:
\exmpl{
Let us first reconsider the bouncing ball example with uncertainty in the coefficient of restitution - which enters into the reset map. Assume the parameter lies on an interval $p\in (0,C)$. Then the model is give by the tuple:
\[H_B(p) = (Q,E,D,F,G,R(p))\]
where
\begin{itemize}
\item  $Q = \{1\}$, which provides the discrete state
\item  $E = \{(1,1)\}$, which is the single edge from $q_0$ to itself
\item  $D := \{x\in\Realsn{2}: x_1\geq 0\}$ provides the domain. Thus, $g_{q_0} = x_1$.
\item $F = f(x)$ provides a vector field mapping $D$ to itself, and where
\[\dot{x} = f(x) = \rmatrix{c}{x_2\\-g}\]
\item $G = \{x\in\Realsn{2}: x_1=0,\;x_2\leq 0\}$ provides the guard. Thus, $h_{(1,1),0} = x_1$, and $h_{(q_0,q_0),1} = -x_2$.
\item $R(p) = \phi(x) = [0,-px_2]^T$ provides the reset map.
\end{itemize}
}\par
{\bf Results:} We would like to prove stability of the Zeno equilibrium for all $p\in [0,C]$ for some $C$. To do this we define the polynomial $\pi(p)_{11}=-p(p-1)$ which yields the uncertainty set  $P=\{p\in\Reals: \tilde{p}(p) := p(p-C)\leq 0\}=[0,C]$. As before, the Zeno equilibrium is $z=[0, \, 0]^T$ and we choose $W_1:=\{x\in D_1 \, : \, w(x) = 25-x_1^2 \geq 0\}$.
From the previous example, we expect that this Zeno equilibrium is stable for $C<1$. Using a 4th degree polynomials for $V(x)$ and the SOS and polynomial multipliers, we performed a bisection search for the maximum $C$ for which this parameterized hybrid model is stable. Our experiments were able to verify Zeno stability for up to $C=0.99$ - which agrees well with the known analytical result.

\exmpl{
Next, we consider a hybrid model with uncertainty in the switching surface - which determines the domains and guard set. Specifically, consider the vector-field in figure~\ref{fig:zenoswitch1}. In this example, the lower switching surface is fixed and the upper surface is allowed to vary between $0^\circ$ and $90^\circ$. The uncertainty is parameterized by $p\in [0,\infty)$ which represents the slope of the upper switching surface. This is described by the parameterized hybrid model $H(p)=(Q,E,D(p),F,G(p),R)$ where
\begin{itemize}
\item $Q = \{1,2\}$
\item $E = \{(1,2),(2,1)\}$
\item $D(p) = \{D_1(p),D_2(p)\}$ where
\aligneq{
&D_1(p) := \{x\in \Realsn{2}: x_1+x_2\geq 0,\, px_1-x_2\geq 0\}\\
&D_2(p) := \{x\in \Realsn{2}: x_1+x_2\leq 0\;\; \text{OR}\;\; px_1-x_2\leq 0\}
}
\item $F = \{f_1,f_2\}$ where
\aligneq{
&f_1 = \rmatrix{c}{-0.1\\
2}\\
&f_2 = \rmatrix{c}{- x_2 - x_1^3\\
x_1}}
\item $G(p) = \{G_{12},G_{21}(p)\}$ where
\aligneq{
&G_{12}(p)=x_2-px_1=0\\
&G_{21} := \{x\in \Realsn{2}: x_1+x_2 = 0\}
}
\item $R= \{\phi_{12}(x),\phi_{21}(x)\}$ where each $\phi_{ij}(x) = x$.
\end{itemize}

\begin{figure}[h]
\centering
\includegraphics[width=0.8\columnwidth]{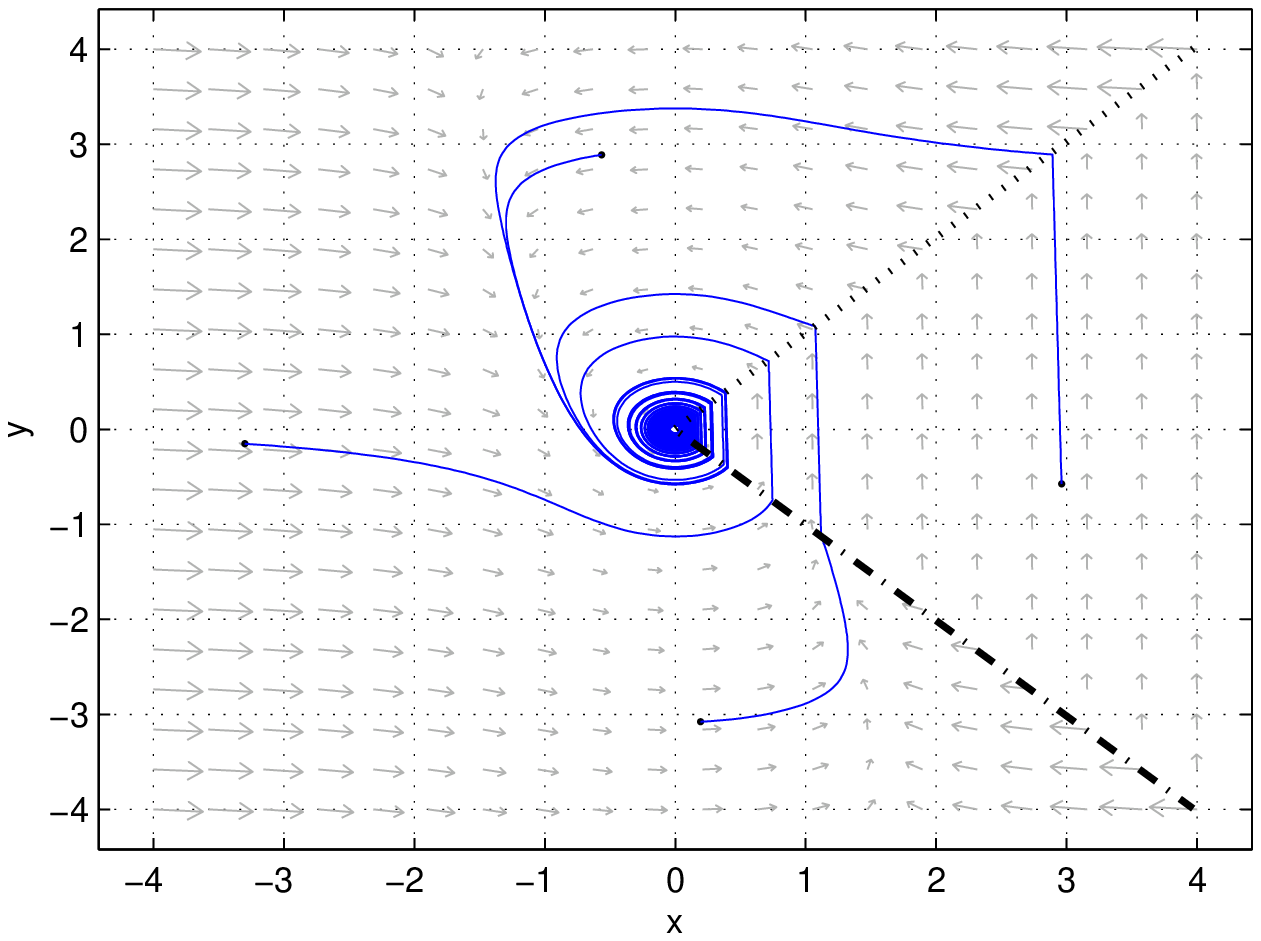}
\caption{Trajectories of Hybrid System in Example 5 with p=1. Dotted line indicates $G_{12}$ and dash-dotted line indicates $G_{21}$}

\label{fig:zenoswitch1}
\end{figure}

\begin{figure}[h]
\centering
\includegraphics[width=0.8\columnwidth]{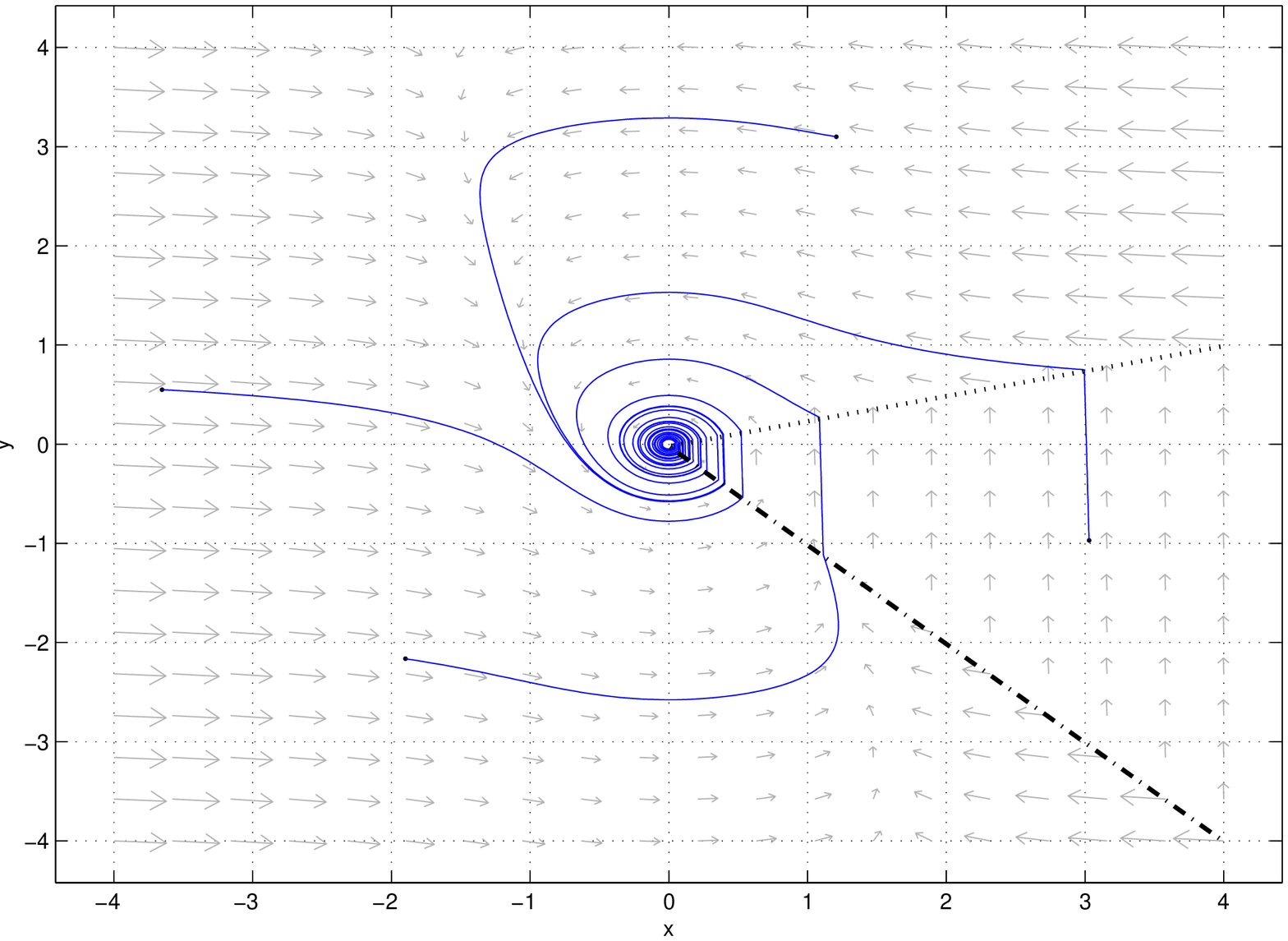}
\caption{Trajectories of Hybrid System in Example 5 with p=4. Dotted line indicates $G_{12}$ and dash-dotted line indicates $G_{21}$}

\label{fig:zenoswitch3}
\end{figure}

\begin{figure}[h]
\centering
\includegraphics[width=0.8\columnwidth]{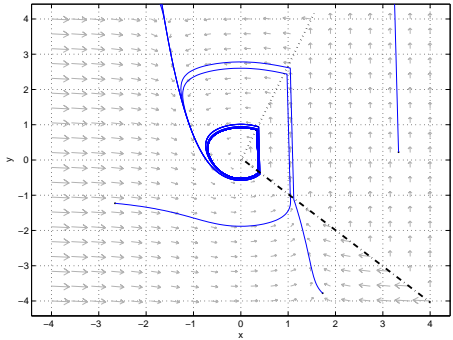}
\caption{Trajectories of Hybrid System in Example 5 with p=0.4. Dotted line indicates $G_{12}$ and dash-dotted line indicates $G_{21}$}

\label{fig:zenoswitch2}
\end{figure}

{\bf Results:} In this example, we use $z=\left\{\bmat{0&0}^T,\bmat{0&0}^T\right\}$ and $W_q(p):=\{x\in D_q(p) \, : \, w(x) = 25-x_1^2 - x_2^2 \geq 0\}$. Simulation indicates the origin is Zeno stable for $p>1$. For $p\in(-0.1,1)$, trajectories converge to a stable limit cycle, as is illustrated in Figure~\ref{fig:zenoswitch2}. If $p\leq -0.1$, then we find the system will no longer stable in any sense.

The first difficulty with this example is that the domain, $D_2(p)$ is NOT a semialgebraic set. To resolve this problem, we represent $D_2$ as the union of two semialgebraic sets $D_{21}$ and $D_{22}$,where
\aligneq{
&D_{21} := \{x\in \Realsn{2}: -px_1+x_2\geq 0\}\\
&D_{22} := \{x\in \Realsn{2}: x_1+x_2\leq 0\}.
}
In this case, for $q=2$, Conditions~\ref{eq:NC1r},~\ref{eq:NC3r} and~\ref{eq:NC4r} must be enforced separately for $W_2(p):=\{x\in D_{21}(p) \, : \, w(x) = 25-x_1^2 - x_2^2 \geq 0\}$ and $W_2(p):=\{x\in D_{22}(p) \, : \, w(x) = 25-x_1^2 - x_2^2 \geq 0\}$. Practically, this means that we have three additional constraints in Feasibility problem 2 corresponding to Constraints~\ref{eq:R1},~\ref{eq:R2} and~\ref{eq:R3} applied to both $g_{21}(x,p)=-x_1-x_2$ and $g_{21}(x,p)=-px_1+x_2$.

We represent the set of uncertainties as
$P:=[C,0]=\{p\in\Reals: \pi(p):=p-C>0\}$, where we must specify the lower $C$. The goal is to find the smallest $C$ such we can prove Zeno stability of $H(p)$ for all $p \in [C,\infty]$. Note that this is actually somewhat challenging as the simplified Positivstellensatz results we discussed earlier only apply to bounded sets. For a fixed polynomial degree, we determine the lowest stable value of $C$ by bisection. As we increase the degree of the polynomials, our lower bound on $C$ improves, as illustrated in Table~\ref{tab:tab1}. Note that we were unable to find a feasible $V_1$ and $V_2$ of degree less than 8 and we were unable to search for polynomials of degree greater than 12 owing to limited computational power.

\begin{table}[h]
\centering
\begin{tabular}{|c|c|}\hline
{\bf Degree of $V_1,V_2$} &{\bf Lower bound on $C$}\\\hline
8 & 2.11\\
10 & 1.87 \\
12 & 1.73\\\hline
\end{tabular}
\caption{Lower bound on $C$ for which $z$ is Zeno stable obtained for different degrees of $V_1$ and $V_2$}
\label{tab:tab1}

\end{table}
}
\vspace{-.1in}
\section{Conclusions}
\vspace{-.1in}
In this paper, we have presented an approach to testing stability of Zeno equilibria for a general class of nonlinear hybrid systems. Our approach is based on application of sum-of-squares optimization to construct high-degree polynomials which satisfy a new class of Lyapunov conditions. This approach can potentially be used to verify convergence on compact sets and accurately estimate domains of attraction. We also consider a class of hybrid systems with parametric uncertainty in the vector field, domain, guard set and reset map and show how our conditions can be applied to these parameterized systems with a semialgebraic uncertainty set. To illustrate this work, we use a number of examples including a parameterized bouncing ball, a variable structure control system, and a Gain-scheduled system, among others. We use our approach to numerically examine the robustness of these Zeno equilibria to uncertainties in the domain, reset map and guard set (switching surface). Our numerical tests indicate convergence of the accuracy of the proposed conditions to the analytic limit.
\vspace{-.1in}
\bibliographystyle{elsarticle-num}
\bibliography{chai_zeno_scl}
\end{document}